\documentclass[11pt,a4paper,equation]{article}
\usepackage[greek,american]{babel}
\usepackage[iso-8859-7]{inputenc}
\usepackage{amssymb,amsfonts}
\usepackage[dvips]{graphicx}
\usepackage{amsmath}
\usepackage{epsfig}
\usepackage{latexsym}
\usepackage{makeidx}
\usepackage{pst-math,pst-xkey}
\newtheorem{lemma}{Lemma}[section]
\newtheorem{theorem}{Theorem}[section]

\numberwithin{equation}{section}
\makeindex
\begin{document}
\date{}
\author{ Aristides V. Doumas$^{1}$ and Vassilis G. Papanicolaou$^{2}$ \\
Department of Mathematics\\
National Technical University of Athens\\
Zografou Campus\\
157 80 Athens, GREECE\\
$^{1}$aris.doumas@hotmail.com \quad $^{2}$papanico@math.ntua.gr}
\title{The logarithmic Zipf version of the coupon collector's problem}
\maketitle
\begin{abstract}
A collector wishes to collect $m$ complete sets of $N$ distinct coupons. The draws from the population are considered to be independent and identical distributed with replacement, and the probability that a type-$j$ coupon is drawn is noted as $p_{j}$. Let $T_{m}(N)$ the number of trials needed for this problem. We present the asymptotics for the expectation (five terms plus an error), the second rising moment (six terms plus an error), and the variance of $T_{m}(N)$ (leading term), as well as its limit distribution as $N\rightarrow \infty$, when 
\begin{equation*}
p_{j}=\frac{a_{j}}{\sum_{j=2}^{N+1} a_{j}}, \,\,\,\text{where}\,\,\, a_{j}=\left(\ln j\right)^{-p}, \,\,p>0. 
\end{equation*}
These ``log-Zipf" classes of coupon probabilities are not covered by the existing literature and the present paper comes to fill this gap. Therefore, we enlarge the classes for which the collector's problem is solved (moments, variance, distribution). 
\end{abstract}

\textbf{Keywords.} Urn problems; coupon collector's problem; double Dixie cup problem; Gumbel
distribution; Laplace method for integrals - Determination of higher order terms; Generalized Zipf law, Eulerian logarithmic integral.\\\\
\textbf{2010 AMS Mathematics Classification.} 60F05; 60F99; 60G70.

\section{Introduction and Motivation}
The coupon collector's problem (CCP) is a classic urn problem of probability theory. It refers to a population whose members are of $N$ different \emph{types} (e.g., baseball cards, viruses, fish, words, etc). For
$1 \leq j \leq N$ we denote by $p_j$ the probability that a member of the population is of type $j$, where $p_j > 0$ and $\sum_{j=1}^{N}p_{j}=1$. The members of the population are sampled independently \textit{with replacement} and their types are recorded. Naturally, the main object of study is the number $T(N)$ of trials needed until all $N$ types are detected (at least once). 
The simple case where all $p_{j}$'s are equal has a long history. It began with A. De Moivre at the eighteenth century and later with P.S. Laplace (see \cite{Ho}, \cite{D-H}).\\
In the recent years D.J. Newman and L. Shepp studied the more general problem where the collector's goal is to complete $m$ sets of all $N$ existing different coupons (still uniformly distributed), \cite{N-S}. This problem is known as the double Dixie cup problem due to a successful marketing policy of the Dixie Cup Company, (see \cite{Ma}). Let $T_{m}(N)$ be the number of trials needed for this case. The main result of \cite{N-S} was that for any fixed $m$
\begin{equation}
E\left[\, T_m(N)\,\right]= N \ln N + \left(m-1\right) N \ln \ln N + N C_m + o(N)
\label{1}
\end{equation}
as $N\rightarrow \infty$, where $C_{m}$ is a constant depending on $m$. 
Soonafter, P. Erd\H{o}s and A. R\'{e}nyi went a step further and determined the limit distribution of $T_{m}(N)$, as well as the exact value of the constant $C_{m}$, see \cite{E-R}. They proved that 
\begin{equation}
C_m = \gamma - \ln \left(m-1\right)!,\label{2}
\end{equation}
where $\gamma=0.5772\cdots$ is the Euler-Mascheroni constant, and that for every real $y$ the following limiting result holds:
\begin{equation}
\lim_{N \rightarrow \infty} P\left\{\frac{T_m(N) - N \ln N - (m - 1) N \ln \ln N + N \ln (m-1)!}{N} \leq y\right\}
= e^{-e^{-y}}
\label{333}
\end{equation}
(the right-hand side of \eqref{333} is the standard Gumbel distribution function). 
For the case of unequal coupon probabilities R.K. Brayton (1963) under the quite restrictive assumption of ``nearly equal coupon probabilities", namely
\begin{equation*}
\lambda(N) :=
\frac{\max_{1\leq j\leq N}{\left\{p_{j}\right\}}}{\min_{1\leq j\leq N}{\left\{p_{j}\right\}}}\leq M < \infty, \qquad\text{independently of $N$,}
\end{equation*}
employed  the formulae
\begin{align}
E[T_m(N)]&=\int_{0}^{\infty}\left\{1-\prod_{j=1}^{N}\left[1-S_{m}(p_{j}t) e^{-p_{j}t} \right]\right\} dt,
\label{5}
\\
E\left[T_m(N)\left(T_{m}(N)+1\right)\right]&=2\int_{0}^{\infty}
\left\{1-\prod_{j=1}^{N}\left[1-S_{m}(p_{j}t) e^{-p_{j}t}\right]\right\} t dt
\label{5A}
\end{align}
and obtained \cite{B} detailed asymptotics of the expectation $E[T_m(N)]$ and the second rising moment $E\left[T_m(N)\left(T_{m}(N)+1\right)\right]$. Here and in what follows $S_{m}(y)$ denotes the $m$-th partial sum of $e^{y}$, namely
\begin{equation}
S_{m}(y) := 1+y+\frac{y^{2}}{2!}+\cdots+\frac{y^{m-1}}{\left(m-1\right)!}=\sum_{l=0}^{m-1}\frac{y^l}{l!}\,.
\label{7}
\end{equation}
As for the asymptotics of the variance, he only did the
case $m=1$, where he found the formula
\begin{equation*}
V\left[\, T_1(N)\,\right]
= N^{2}\left[\frac{\pi^{2}}{6} + O\left(\frac{\ln \ln \ln N}{\ln \ln N}\right)\right] \quad \text{as}\quad N\rightarrow \infty.
\end{equation*}
For the case of unequal coupon probabilities and for $m=1$, general results have been published in \cite{DP} and \cite{DPM}, while for general (however fixed) values of $m$ a paper of ours has been recently uploaded in the \textit{arxiv}, \cite{MSETS}. Since our motivation arises from these works we will briefly present their results. Let $\alpha =\{a_{j}\}_{j=1}^{\infty }$ be a sequence of strictly positive
numbers. Then, for each integer $N > 0$, one can create a probability measure
$\pi _N =\{p_1,...,p_N\}$ on the set of types $\{1,...,N\}$ by taking
\begin{equation}
p_j = \frac{a_j}{A_N},
\qquad \text{where}\quad
A_N = \sum_{j=1}^N a_j.
\label{8}
\end{equation}
Notice that $p_j$ depends on $\alpha $ and $N$, thus, given $\alpha $, it
makes sense to consider the asymptotic behavior of $E\left[\, T_m(N)\,\right]$, $E\left[\,T_m(N)\left(T_{m}(N)+1\right)\,\right]$, and $V\left[\, T_m(N)\,\right]$ as $N\rightarrow \infty$. It follows that 
\begin{equation}
E\left[\, T_m(N)\,\right] =A_{N}\,E_{m}(N;\alpha),
\label{12}
\end{equation}
\begin{equation}
E\left[\,T_m(N)\left(T_{m}(N)+1\right)\,\right] =A^{2}_{N}\,Q_{m}(N;\alpha),
\label{15}
\end{equation}
where
\begin{align}
E_{m}(N;\alpha ):&=\int_{0}^{\infty}\left[1-\prod_{j=1}^{N}\bigg(1-e^{-a_{j}t}\,S_{m}\left(a_{j}t\right)\bigg)\right]dt, \label{9} \\
Q_{m}(N;\alpha ):&=2\int_{0}^{\infty}t\left[1-\prod_{j=1}^{N}\bigg(1-e^{-a_{j}t}\,S_{m}\left(a_{j}t\right)\bigg)\right]dt. \label{13}
\end{align}
Let
\begin{equation}
L_{1}(\alpha;m ):=\lim_{N}E_{m}(N;\alpha )\,\,\,\text{and}\,\,\,L_{2}(\alpha;m ):=\lim_{N}Q_{m}(N;\alpha ).
\label{17}
\end{equation}
The sequences $\alpha=\left\{a_{j}\right\}_{j=1}^{\infty}$ were separated as follows:
\begin{equation*}
\text{\textbf{(Case I)}}\,\,\,\sum_{j = 1}^\infty e^{-a_j \tau} < \infty\,\,\,\,\text{for some}\,\, \tau>0.
\end{equation*}
Notice that Case I is equivalent to $L_1 (\alpha;m)<\infty$ and $L_2 (\alpha;m)< \infty$.
As it turned out the leading term of both the expectation and the second (rising) moment of $T_{m}(N)$ is enough to obtain the leading asymptotics of its variance. As for the distribution of $T_m(N)$, for all $s \in [0, \infty)$ one has
\begin{equation*}
P\left\{\frac{T_m(N)}{A_N} \leq s \right\}\rightarrow F(s) := \prod_{j=1}^{\infty}\left[1 - S_{m}(a_j s) e^{-a_j s} \right],
\qquad
N \rightarrow \infty,
\end{equation*}
where $S_m(\,\cdot \,)$ is given by (\ref{7}).

\smallskip 

Examples of sequences falling in this case are $a_{j}=j^{p}$, $p>0$ (for $p=1$ we have the so-called \textit{linear case}), $b_{j}=e^{p j}$, $p>0$, and $c_{j}=j!$.
\begin{equation*}
\text{\textbf{(Case II)}}\,\,\,\sum_{j = 1}^\infty e^{-a_j \tau} = \infty\,\,\,\,\text{for all}\,\, \tau>0,
\end{equation*}
which is equivalent to $L_1 (\alpha;m)=L_2 (\alpha;m)=\infty$. In order to make some progress one has to make some assumptions for the sequence $\alpha=\left\{a_{j}\right\}_{j=1}^{\infty}$. If we write $a_{j}$ as
\begin{equation}
a_{j}=f(j)^{-1},
\label{aj}
\end{equation}
where
\begin{equation*}
f(x) > 0 \qquad \text{and} \qquad f'(x) > 0,
\end{equation*} 
and assume that $f(x)$ possesses three derivatives satisfying the following conditions as $x\rightarrow \infty$:
\begin{align}
\text{(i) }& f(x)\rightarrow \infty,
& &\text{(ii) } \frac{f^{\prime }(x)}{f(x)}\rightarrow 0,\nonumber \\
\text{(iii) } &\frac{f^{\prime \prime}(x)/f^{\prime }(x)}{f'(x)/f(x)} = O\left(1\right),
& &\text{(iv) } \frac{f^{\prime \prime\prime}(x)\;f(x)^{2}}{ f^{\prime }(x)^{3}} = O\left(1\right),
\label{C1}
\end{align}
then, the asymptotics of the expectation of $T_{m}(N)$ (up to the fifth term), and the asymptotics of its second rising moment (up to the sixth term) were obtained. These results were needed for the leading asymptotics of the variance $V[\,T_{m}(N)\,]$ to appear. As for the limiting distribution as it turned out the random variable $T_{m}(N)$ (under the appropriate normalization) converges in distribution to a Gumbel random variable.

\bigskip

\textbf{Remark 1.} Roughly speaking, $f(\cdot)$ belongs to the class of positive and strictly increasing functions, which grow to $\infty$
(as $x \rightarrow \infty$) \textit{slower than exponentials, but faster than powers of logarithms}.

\smallskip

In particular, $(ii)$ is a sub exponential condition. Conditions $(iii)$ (mainly) and $(iv)$ interpret the above remark for the growth of $f(\cdot)$. These conditions are satisfied by a variety of commonly used functions. For example,
\begin{equation*}
f(x) = x^p (\ln x)^q, \quad p > 0,\ q \in \mathbb{R},\qquad \qquad
f(x) = \exp(x^{r}),\quad 0 < r < 1,
\end{equation*}
or various convex  combinations of products of such functions.\\
In particular when
\begin{equation*}
f(x)=x^{p},\,\, p>0
\end{equation*}
that is the coupon probabilities are
\begin{equation}
p_{j}=\frac{a_{j}}{\sum_{j=1}^{N}a_{j}},\,\,\,\,a_{j}=\frac{1}{j^{p}},\,\,\, p>0 \label{Z}
\end{equation}
we have the so-called \textit{generalized Zipf distribution}, a surprising law, which have attracted the interest of many researchers, mainly due to its application in computer science and linguistics (the literature on the Zipf law is extensive). In reference to the CCP the standard Zipf distribution (that is the case where $p=1$) and when $m=1$, the asymptotics of the expectation (leading term) of $T_{1}(N)$, was first studied by Flajolet \textit{et al}, see \cite{F-G-T}.\\
To summarize, we have an answer for the asymptotics of the expectation and the second rising moment of $T_{m}(N)$, as well as the leading asymptotics of the variance $V[\,T_{m}(N)\,]$, and its limiting distribution for rich classes of coupon probabilities. Moreover, even exponential sequences belong to the set of classes of sequences, for which we are able to solve our problem. For example the sequence $\beta = \{e^{-pj} \}_{j=1}^{\infty}$, $p > 0$ falls into Case II; but condition $(ii)$ of (\ref{C1}) is violated. However, if one considers the sequence $\alpha = \{e^{pj} \}_{j=1}^{\infty}$ it is immediate  that $\alpha$ and $\beta$ produce the same coupon probabilities, and since $\alpha$ falls into Case I, a solution to our problem exists.\\
The question arises naturally: can we extend the classes of functions $f(\cdot)$? What happens if our functions grows as powers of logarithms?

\bigskip

\textbf{Problem.} What can be said about the moments, the variance, and the distribution of the random variable $T_{m}(N)$, when $f(x)= \ln x$, or more generally when $f(x)= (\ln x)^{-p}$, $p>0$? In other words what can be said for the case the coupon probabilities satisfy:
\begin{equation}
p_{j}=\frac{a_{j}}{\sum_{j=2}^{N+1} a_{j}}, \,\,\,\text{where}\,\,\, a_{j}=\left(\ln j\right)^{-p}, \,\,p>0. \label{cc}
\end{equation}
\textbf{Remark 2.} Formulae (\ref{Z}) and (\ref{cc}) are explaining the title of this paper.
\section{Discussion and main results}
Consider the case $a_{j}=\left(\ln j\right)^{-p}, \,\,p>0$. Clearly,
\begin{equation*}
\sum_{j=2}^{\infty}e^{-\tau\left(\ln j\right)^{-p}}=\infty \quad \text{for all}\quad \tau>0.
\end{equation*}
Therefore, these sequences fall into Case II. However, conditions $(iii)$ and $(iv)$ of (\ref{C1}) are violated. In view of (\ref{cc}), (\ref{12}), and (\ref{5}) we get
\begin{align}
E[\,T_m(N)\,]&=\left(\sum_{j=2}^{N+1}\left(\ln j\right)^{-p}\right)\int_{0}^{\infty}\left\{1-\prod_{j=2}^{N+1}
\left[1-S_{m}\bigg(t \left(\ln j\right)^{-p}\bigg) e^{-t\, \left(\ln j\right)^{-p}} \right]\right\} dt.
\label{exp}
\end{align}
\textbf{Remark 3.} From here and in what follows we replace $N+1$ by $N$, in both the sum and the integral above, without loss of information regarding the asymptotics of $E[\,T_m(N)\,]$.\\\\
The sum $\sum_{j=2}^{N}\left(\ln j\right)^{-p}$ in (\ref{exp}) is easy to handle. In fact one may easily obtain its full asymptotic expanssion by using the Euler-Maclaurin summation formula, and hence the associated integral $\int_{j=2}^{N}\left(\ln x\right)^{-p}dx$, and then repeated integration by parts, (see \cite{B-O}). In particular, for $p=1$ we get the so-called \textit{offset logarithmic integral} or \textit{Eulerian logarithmic integral}, which is a very good approximation to the number of prime numbers less than $N$ (i.e., $\pi (x)\sim \int_{j=2}^{N}\left(\ln x\right)^{-p}dx$). We get
\begin{equation}
A_N = \sum_{j=2}^N \frac{1}{(\ln j)^p} = \frac{N}{(\ln N)^p} + \frac{p\,N}{(\ln N)^{p+1}}
+\frac{p\left(p+1\right)\,N}{(\ln N)^{p+2}} + O\left( \frac{N}{(\ln N)^{p+3}} \right).
\label{SD4}
\end{equation}
The integral appearing in (\ref{exp}) is $E_{m}(N;\alpha )$ of (\ref{13}) and is our main task. Our approach lies in three steps.\\
\underline{Step 1} is a change of variables 
\begin{equation*}
t=g(N)\,s\end{equation*}
where 
\begin{equation*}
\lim_{N} g(N)=\infty.
\end{equation*}
There are maybe infinite choices for $g(N)$, but a convinient one is 
\begin{equation*}
g(N)=(\ln N)^{p+1},
\end{equation*}
which makes things simpler by invoking (\ref{SD4}). Thus,
\begin{align}
E[\,T_m(N)\,]&=\bigg(N\ln N + p\,N+p\left(p+1\right)\frac{N}{\ln N}+ O\left( \frac{N}{\left(\ln N\right)^{2}}\right)\bigg) \nonumber\\ 
&\times \int_{0}^{\infty}\left\{1-\exp \Bigg( \sum_{j=2}^{N}\ln
\left[1-S_{m}\bigg(\frac{(\ln N)^{p+1}}{ \left(\ln j\right)^{p}}s\bigg) \exp{\bigg(-\frac{(\ln N)^{p+1}}{ \left(\ln j\right)^{p}}s\bigg)} \right]\Bigg)\right\} ds.
\label{EXP}
\end{align}
\underline{Step 2}. The asymptotics (as $N\rightarrow \infty$) of the integral
\begin{equation}
I_{k}(N):=\int_{2}^{N} \exp{\bigg(-\frac{(\ln N)^{p+1}}{ \left(\ln x\right)^{p}}\,s\bigg)}
\frac{dx}{\left(\ln x\right)^{kp}}, \quad k=0,1,\cdots,m-1,\,\,\,p>0.
  \label{LL2}
\end{equation}
\begin{lemma}
\begin{equation*}
I_{k}(N)=N^{1-s}\left(\ln N\right)^{-k p}\left[\frac{1}{1+ps}+\frac{k p}{\left(1+ps\right)^{2}\ln N}-\frac{p\left(p+1\right)s}{\left(1+ps\right)^{3}\ln N}\left(1+O\left(\frac{1}{\ln N}\right)\right)\right],
\end{equation*}
uniformly in $s\in [s_{0},\infty)$, for any fixed $s_{0} > 0 $.
\end{lemma}
All the proofs of this paper are gathered in Section 3. For now we only wish to note that the main tool to estimate the integral above is the Laplace method for integrals for the determination of higher order terms. Hence,
\begin{equation}
\lim_{N}\int_{2}^{N}\exp{\bigg(-\frac{(\ln N)^{p+1}}{ \left(\ln x\right)^{p}}s\bigg)} S_{1}\bigg(\frac{(\ln N)^{p+1}}{ \left(\ln j\right)^{p}}s\bigg)dx =\left\{
\begin{array}{rcc}
\infty,&  \text{if } s<1, \\
(1+p)^{-1},& \text{if } s=1, \\
0,&   \text{if } s> 1,
\end{array}
\right. \label{LL2a}
\end{equation}
while for $m\geq 2$
\begin{equation}
\lim_{N}\int_{2}^{N}\exp{\bigg(-\frac{(\ln N)^{p+1}}{ \left(\ln x\right)^{p}}s\bigg)} S_{m}\bigg(\frac{(\ln N)^{p+1}}{ \left(\ln j\right)^{p}}s\bigg)dx  
=\left\{
\begin{array}{rcc}
\infty,&  \text{if }& s\leq 1, \\
0,&   \text{if }& s> 1,
\end{array}
\right. 
\label{LL1}
\end{equation}
Now from the comparison of sums and integrals it follows that the limits above are valid, if the integral is replaced by the associated sum. Moreover, from the Taylor expansion for the logarithm, namely $\ln(1-x)\sim-x$ as $x\rightarrow0$, one gets the corresponding  limits, e.g. for all $m\geq 2$
\begin{equation}
\lim_{N}\sum_{j=2}^{N}\ln \left[1-S_{m}\bigg(\frac{(\ln N)^{p+1}}{ \left(\ln j\right)^{p}}s\bigg) \exp{\bigg(-\frac{(\ln N)^{p+1}}{ \left(\ln j\right)^{p}}s\bigg)} \right] =\left\{
\begin{array}{rc}
-\infty,&  \text{if } s<1 \\
0,&   \text{if } s\geq1.
\end{array}
\right.  \label{SL1A}
\end{equation}
The limit above drives us to \underline{Step 3}. This is actually a method we proposed recently in \cite{DP}. We do not claim that this method is new, but even though there is \textit{no guarantee} that it can be applied in our problem (since conditions (\ref{C1}) are violated), it turns out that it is leads to a solution. We will briefly discuss it here and complete the proof in the next section. Let us denote by $\tilde{E}_{m}(N;\alpha)$ the integral appearing in (\ref{EXP}). For any given $\varepsilon \in (0,1)$ one has
\begin{align}
\tilde{E}_{m}(N;\alpha)= \left[\,1+\varepsilon -I_1 (N)-I_2 (N)+I_3 (N)\,\right],
\label{b3a}
\end{align}
where
\begin{align}
I_1 (N):&= \int_0^{1-\varepsilon} e^{M_{m}(N;s)}\, ds,\label{I1}\\
I_{2}(N): &= \int_{1-\varepsilon}^{1+\varepsilon } e^{M_{m}(N;s)}\, ds,\label{I2}\\
I_{3}(N): &=\int_{1+\varepsilon}^{\infty } 1-e^{M_{m}(N;s)}\, ds,\label{I3}
\end{align}
and
\begin{equation}
M_{m}(N;s) :=\sum_{j=2}^N \ln \left[1-S_{m}\bigg(\frac{(\ln N)^{p+1}}{ \left(\ln j\right)^{p}}s\bigg) \exp{\bigg(-\frac{(\ln N)^{p+1}}{ \left(\ln j\right)^{p}}s\bigg)} \right].
\label{AsN}
\end{equation}
The heart of Step 3 is that $I_{3}(N)$ and $I_{1}(N)$ are dominated by the sixth term in the asymptotics of $I_{2}(N)$ as $N\rightarrow \infty$. Intuitively one expects that the main contribution of $\tilde{E}_{m}(N;\alpha)$ should come from $I_{2}(N)$ (due to the limit of (\ref{SL1A}), but it turns out that $I_{2}(N)$ is much more important. The analysis of $I_{2}(N)$ lies in Lemma 2.1(critical contribution), as well as in classical techniques of asymptotic analysis. The computations needed are often, quite involved.
\begin{theorem} \textbf{(Main result I)}\\
Let $T_m(N)$ the number of trials a collector needs to complete $m$ sets of $N$ different types of coupons with replacement. If the coupon probabilities satisfy
\begin{equation*}
p_{j}=\frac{a_{j}}{\sum_{j=2}^{N}a_{j}}, \quad\text{where}\quad a_{j}=\left(\ln j\right)^{-p}, \,\,\,p>0
\end{equation*} 
then, the asymptotics of the average of $T_m(N)$ (as $N\rightarrow \infty$) satisfy
\begin{align}
E\left[\, T_m(N)\,\right] = N \ln N &+ \left(m-1\right)N \ln\ln N + \left[\,p+\gamma-\ln\left(m-1\right)!-\ln\left(p+1\right)\,\right]\,N \nonumber \\
&-(m-1)\left[\frac{p}{p+1}-\left(m-1\right)-p\right]\,\frac{\ln\ln N}{\ln N}\,N\nonumber \\
&+N\left[p\left(p+1\right)-p\,\bigg(\ln\left(m-1\right)!+\ln\left(p+1\right)-\gamma \bigg)\right.\nonumber\\
&\left.\,\,\,\,\,\,\,\,\,\,\,\,\,\,\,\,\,\,\,-\bigg(\frac{p}{p+1}-\left(m-1\right)\bigg)\times\left[\gamma-\ln\left(m-1\right)!-\ln\left(p+1\right)\right.\right.\nonumber\\
&\left.\left.\,\,\,\,\,\,\,\,\,\,\,\,\,\,\,\,\,\,\,\,\,\,\,\,\,\,\,\,\,\,\,\,\,\,\,\,\,\,-\frac{1}{\left(p+1\right)^{2}}\left(\frac{m-1}{p+1}-\frac{p+1}{p}-3\left(\frac{p}{p+1}\right)^{2}\right)\right]\right]\nonumber\\
&+O\left(\frac{\ln\ln N}{\left(\ln N\right)^{2}}\,N\right),
\label{R1}
\end{align}
where $\gamma$ is, as usual, the Euler-Mascheroni constant.
\end{theorem}   
\textbf{Remark 4.} Notice that the expected value in \eqref{R1} is slightly bigger than the corresponding expected value for the case of equal coupon probabilities
(recall \eqref{1}--\eqref{2}), due to the term $p - \ln(p+1)$ which is strictly positive for all $p > 0$. This is in accordance with the statement: For fixed positive integers $m$ and $N$, the case of equal probabilities, has the property that it is the one with the stochastically smallest $T_m(N)$. This result is due to \cite{MO}.
\begin{theorem} \textbf{(Main result II)}\\
For the second (rising) moment of the random variable $T_m(N)$ we have the following asymptotic expression as $N\rightarrow \infty$
\begin{align}
E\left[\,T_m(N)\left(T_{m}(N)+1\right)\,\right]& = N^{2} \left(\ln N\right)^{2} +2\left(m-1\right)N^{2}{\ln N}\left(\ln\ln N\right)\nonumber \\
&+2 \left[\,p+\gamma-\ln\left(m-1\right)!-\ln\left(p+1\right)\,\right]\,N^{2}\ln N \nonumber \\
&+\left(m-1\right)^{2}\,N^{2}\left(\ln\ln N\right)^{2}\nonumber \\
&-2\left(m-1\right)\left(\frac{p}{p+1}-\left(m-1\right)-\gamma-2p\right.\nonumber\\
&\left.\,\,\,\,\,\,\,\,\,\,\,\,\,\,\,\,\,\,\,\,\,\,\,\,\,\,\,\,\,\,\,\,\,\,\,\,\,\,\,\,\,\,\,\,\,\,\,\,\,\,\,\,+\ln\left(m-1\right)!+\ln\left(p+1\right)\right)N^{2}\ln\ln N\nonumber \\
&+N^{2}\left[p^{2}+2p\left(p+1\right)-2\left(2p+\gamma\right)\,\bigg(\ln\left(m-1\right)!+\ln\left(p+1\right)\bigg)\right.\nonumber\\
&\left.\,\,\,\,\,\,\,\,\,\,\,\,\,\,\,\,\,\,\,\,\,\,\,\,+4p\gamma-\bigg(\ln\left(m-1\right)!+\ln\left(p+1\right)\bigg)^{2}+\gamma^{2}+\frac{\pi^{2}}{6}\right.\nonumber\\
&\left.\,\,\,\,\,\,\,\,\,\,\,\,\,\,\,\,\,\,\,\,\,\,\,\,-2\bigg(\frac{p}{p+1}-\left(m-1\right)\bigg)\times\left[\gamma-\ln\left(m-1\right)!-\ln\left(p+1\right)\right.\right.\nonumber\\
&\left.\left.\,\,\,\,\,\,\,\,\,\,\,\,\,\,\,\,\,\,\,\,\,\,\,\,\,\,\,\,\,\,\,\,\,\,\,\,\,\,\,\,\,\,-\frac{1}{\left(p+1\right)^{2}}\left(\frac{m-1}{p+1}-\frac{p+1}{p}-3\left(\frac{p}{p+1}\right)^{2}\right)\right]\right]\nonumber\\
&+O\left(\frac{\left(\ln\ln N\right)^{2}}{\ln N}\,N^{2}\right).
\label{R2}
\end{align}
\end{theorem}
\begin{theorem} \textbf{(Main result III)}\\ 
Let $T_m(N)$ the number of trials a collector needs to complete $m$ sets of $N$ different types of coupons with replacement ($m$ is a fixed positive integer). When the coupon probabilities satisfy
\begin{equation*}
p_{j}=\frac{a_{j}}{\sum_{j=2}^{N}a_{j}}, \quad\text{where}\quad a_{j}=\left(\ln j\right)^{-p}, \,\,\,p>0
\end{equation*} 
we have as $N\rightarrow \infty$
\begin{equation}
V\left[\,T_m(N)\,\right] \sim \frac{\pi^2}{6}\;N^2
\label{FINAL}
\end{equation}
{independently} of the value of the positive integer $m$.
\end{theorem}
Having detailed asymptotics for $E\left[\,T_m(N)\,\right]$ and the leading asymptotics for the variance $V\left[\,T_m(N)\,\right]$ we take advantage of a well known but very general limit theorem of P.Neal (see Section 3), and present the following
\begin{theorem} \textbf{(Main result IV)}\\
Suppose the coupon probabilities $p_{j}$ come from the sequence $\alpha = \{a_j = (\ln j)^{-p}\}_{j=2}^{\infty}$ for some $p > 0$, $p_{j}=a_{j}/\sum_{j=2}^{N}a_{j}$. Then, for all $y \in \mathbb{R}$ and for all positive integer $m$ we have as $N \rightarrow \infty$
\begin{equation}
P\left\{\frac{T_m(N) - N \ln N - (m-1) N \ln\ln N - \left[\gamma + p - \ln\bigg((p+1)(m-1)!\bigg) \right] N }{N} \leq y \right\}\rightarrow e^{-e^{-y}}.
\label{SD1a}
\end{equation}
That is, the random variable $T_{m}(N)$ (under the normalization above) converges in distribution to a Gumbel random variable.
\end{theorem}
\subsection{Final comments}
The main task of this paper is to enlarge the classes of coupon probabilities for which we have an answer to the collector's problem (and in general for the Dixie cup problem) for the average, the variance and the limiting distribution. Since the full asymptotic expansion of $\sum_{j=2}^{N}\left(\ln j\right)^{-p}$ is available our approach is analytic (continuous). We approximate sums by integrals. For example, a key formula is (\ref{AsN}), which is valid for $m\geq2$:
\begin{equation*}
\lim_{N}\sum_{j=2}^{N}\ln \left[1-S_{m}\bigg(\frac{(\ln N)^{p+1}}{ \left(\ln j\right)^{p}}s\bigg) \exp{\bigg(-\frac{(\ln N)^{p+1}}{ \left(\ln j\right)^{p}}s\bigg)} \right] =\left\{
\begin{array}{rc}
-\infty,&  \text{if } s<1 \\
0,&   \text{if } s\geq1.
\end{array}
\right.  
\end{equation*}
As for the corresponding integrals we apply the Laplace method for the determination of higher order terms. The analysis of these integrals is complicated. We build on the method proposed in previous works of ours even though the original conditions are violated and one would expect that this approach does not guarantee a path to a solution. We believe that this method could be valuable for future researchers in order to further enlarge the classes of distributions for this problem.\\
Let us now comment on the moments of the random variable $T_{m}(N)$. In view of (\ref{I2}) and (\ref{I5}) (see Section 3), the key integral for the $r$ rising moment of $T_{m}(N)$ should be
\begin{align*}
I(N): &= \int_{1-\varepsilon}^{1+\varepsilon } s^{r-1}\,e^{M_{m}(N;s)}\, ds,
\end{align*}
where 
\begin{equation*}
M_{m}(N;s) :=\sum_{j=2}^N \ln \left[1-S_{m}\bigg(\frac{(\ln N)^{p+1}}{ \left(\ln j\right)^{p}}s\bigg) \exp{\bigg(-\frac{(\ln N)^{p+1}}{ \left(\ln j\right)^{p}}s\bigg)} \right].
\end{equation*}
To give closure let us illustrate a concrete instance of Main result IV motivated by the following example from Feller, \cite{F} (which is also in Durrett, \cite{D}):

\smallskip

\textbf{Example.} What is the probability that in a village of $2190\, (=6 \cdot 365)$ people all birthdays are presented? Is the answer much different for $1825\, (=5\cdot 365)$ people?

\smallskip
We will answer for both cases of uniform and log-Zipf distributions.\\
In the case of equal probabilities we apply the result of P. Erd\H{o}s and A. R\'{e}nyi, see (\ref{333}) and get (since $m=1$) 
\begin{eqnarray*}
P\left(T_{\text{equal}}(365)\leq 2190\right)&=&P\left(\left(T_{\text{equal}}(365)-2153\right)/365\leq37/365\right)  \\
&\approx&\exp(-e^{-0.1014})=\exp(-0.9036)=0.4051.
\end{eqnarray*}
On the other hand
\begin{eqnarray*}
P\left(T_{\text{equal}}(365)\leq 1825\right)&=&P\left(\left(T_{\text{equal}}(365)-2153\right)/365\leq -328/365\right)  \\
&\approx&\exp(-e^{0.8986})=\exp(-2.4562)=0.085.
\end{eqnarray*}
For the case
\begin{equation*}
p_{j}=\frac{a_{j}}{\sum_{j=2}^{366}a_{j}}, \quad\text{where}\quad a_{j}=\left(\ln j\right)^{-1}, 
\end{equation*} 
we apply Main result IV. We have $N=365$, and $N\ln N=2153$,  $(\gamma+1)365=575.684$ and get
\begin{eqnarray*}
P\left(T_{\text{Log Zipf}}(365)\leq 2190\right)&=&P\left(\left(T_{\text{equal}}(365)-2153-575.684\right)/365\leq(-538.684)/365\right)  \\
&\approx&\exp(-e^{1.47585})=0.0126
\end{eqnarray*}
and 
\begin{eqnarray*}
P\left(T_{\text{Log Zipf}}(365)\leq 1825\right)&=&P\left(\left(T_{\text{equal}}(365)-2153-575.684\right)/365\leq(-903.684)/365\right)  \\
&\approx&\exp(-e^{2.47585})=6.84652\times 10^{-6}=0.00000684652
\end{eqnarray*}
Notice that for the equal case, we have the following ratio 
\begin{equation*}
P\left(T_{\text{equal}}(365)\leq 2190\right)/P\left(T_{\text{equal}}(365)\leq 1825\right)=4.77
\end{equation*}
while, for the Log Zipf case we get
\begin{equation*}
P\left(T_{\text{Log Zipf}}(365)\leq 2190\right)/P\left(T_{\text{Log Zipf}}(365)\leq 1825\right)=1840.
\end{equation*}

\section{Proofs}
\textbf{Proof of Lemma 2.1}. From (\ref{LL2}) we easily have
\begin{equation*}
I_{k}(N)=\int_{\ln 2}^{\ln N} \exp{\bigg(-\frac{(\ln N)^{p+1}}{ y^{p}}\,s\bigg)}
\frac{e^{y}}{y^{kp}}dy.
\end{equation*}
The substitution $y=\left(s^{1/p+1} \ln N\right) t $ yields
\begin{equation*}
I_{k}(N)=\frac{s^{\frac{1-k p}{p+1}}}{\left(\ln N\right)^{k p}}\int_{a\, s^{-1/p+1}}^{s^{-1/p+1}} 
\exp{\bigg(s^{1/p+1}\ln N \left(t-t^{-p}\right)\bigg)}
\frac{dt}{t^{kp}}.
\end{equation*}
where $a=\ln2/\ln N$. For convenience we set the integral above as $\tilde{I}_{k}(N)$. Now as long as $s\geq s_{0}>0$ for any fixed $s_{0}$, we have
\begin{equation*}
\lim_{N} s^{1/p+1}\ln N=\infty,\,\,\, \text{for all}\,\, p>0.
\end{equation*}
Moreover, the function 
\begin{equation*}
\phi (t):=t-t^{-p}
\end{equation*}
attains its maximum value at $t_{0}=s^{-1/p+1}$. Hence, only the immediate neighborhood of $t_{0}$ contributes to the full asymptotic expansion of $\tilde{I}_{k}(N)$. Set $h(t):=t^{-k p}$. 
Careful application of Laplace's method for integrals (for the determination of higher-order terms) drives us to approximate $\phi(t)$ by $\phi(t_{0})+(t-t_{0})\phi^{\prime}(t_{0})+\frac{1}{2} (t-t_{0})^{2}\phi^{\prime\prime}(t_{0})$ and $h(t)$ by $h(t_{0})+(t-t_{0})h^{\prime}(t_{0})+\frac{1}{2} (t-t_{0})^{2}h^{\prime\prime}(t_{0})$. Then,
\begin{align*}
\tilde{I}_{k}(N)\sim&\int_{t_{0}-\epsilon}^{t_{0}}  \left[h(t_{0})+(t-t_{0})h^{\prime}(t_{0})+\frac{1}{2} (t-t_{0})^{2}h^{\prime\prime}(t_{0})\right]\\
&\times\exp{\bigg(s^{1/p+1}\ln N \left[\phi(t_{0})+(t-t_{0})\phi^{\prime}(t_{0})+\frac{1}{2} (t-t_{0})^{2}\phi^{\prime\prime}(t_{0})\right]\bigg)}
dt.
\end{align*}
Because $\epsilon$ may be chosen small, we Taylor expand the term 
\begin{equation*}
\exp\left[s^{1/p+1}\ln N \frac{1}{2} (t-t_{0})^{2}\phi^{\prime\prime}(t_{0})\right].
\end{equation*}
Substituting this expansion in the above, then collecting powers of $(t-t_{0})$, and finally, extending the range of integration to $(-\infty,t_{0}]$, yields
\begin{align*}
\tilde{I}_{k}(N)\sim\,&e^{s^{1/p+1}\ln N\,\phi(t_{0})}\int_{-\infty}^{t_{0}}e^{s^{1/p+1}\ln N\,(t-t_{0})\,\phi^{\prime}(t_{0})}\\
&\times\left[h(t_{0})+(t-t_{0})h^{\prime}(t_{0})+\frac{1}{2} (t-t_{0})^{2}\left(h^{\prime\prime}(t_{0})+s^{1/p+1}\ln N\, h(t_{0})\,\phi^{\prime\prime}(t_{0})
\right)+\cdots\right] dt.
\end{align*}
and the proof completes the evaluation of the above integral. For more details on this method, see e.g., \cite{B-O}.\\\\
\textbf{Proof of main result I}.
To analyse (\ref{b3a}) we will start from $I_{2}(N)$ (see (\ref{I2})) and obtain the five first terms in its asymptotic expansion (plus an error). Then we will calculate the leading term of $I_{3}(N)$ and prove that is negligible compared to the sixth term of $I_{2}(N)$ as $N\rightarrow \infty$. Finally, we will estimate the leading term of $I_{1}(N)$, for which we will see that is negligible compared to the leading term of $I_{3}(N)$.\\
Since $\ln (1-x) = -x + O(x^{2})$ as $x\rightarrow 0$, it follows from (\ref{AsN}) and (\ref{7}) that
\begin{align}
M_{m}(N;s)=&-\sum_{k=0}^{m-1}\frac{\left(\ln N\right)^{k(p+1)}\,s^{k}}{k!}
\left(\sum_{j=2}^{N}\left(\ln j\right)^{-k p}\,\exp{-\frac{\left(\ln N\right)^{p+1}}{\left(\ln j\right)^{p}}s}\right) \nonumber \\
&+\sum_{j=1}^{N}O\left(e^{-\frac{2\left(\ln N\right)^{p+1}}{\left(\ln j\right)^{p}}s}\left[S_{m}\left(\frac{(\ln N)^{p+1}}{ \left(\ln j\right)^{p}}s\right)\right]^{2}\right).\label{40}
\end{align}
From the comparison of sums and integrals and Lemma 2.1 (remember that we are interested in $I_{2}(N)$, $s$ is strictly positive and hence we are able to apply Lemma 2.1)
\begin{align}
M_{m}(N;s)=-N^{1-s}\,\sum_{k=0}^{m-1}\frac{\left(\ln N\right)^{k}\,s^{k}}{k!}&\left[\frac{1}{1+ps}+\frac{k p}{\left(1+ps\right)^{2}\ln N}\right.\nonumber\\
&\left.\,\,\,\,\,\,\,\,\,\,\,\,\,\,\,\,\,\,\,\,\,\,-\frac{p\left(p+1\right)s}{\left(1+ps\right)^{3}\ln N}\left(1+O\left(\frac{1}{\ln N}\right)\right)\right]. \label{qw}
\end{align}
Next, we substitute (\ref{qw}) into (\ref{I2})) and apply the change of variables $s=1-t$. Thus,
\begin{align*}
I_{2}(N)=\int_{-\varepsilon}^{\varepsilon}\exp&\left\{-N^{t}\left(\ln N\right)^{m-1}\left(1-t\right)^{m-1}\frac{\left(1-b\right)}{\left(m-1\right)!}\sum_{n=0}^{\infty}
\left(b\,t\right)^{n}+\left(\ln N\right)^{m-2}\frac{\left(1-t\right)^{m-2}}{\left(m-1\right)!}\right.\\
&\left.\times\left[\left(m-1\right)\left(1-b\right)\sum_{n=0}^{\infty}\left(b\,t\right)^{n}+\left(m-1\right)\left(1-b\right)\left(1-t\right)\sum_{n=1}^{\infty} n b^{n}t^{n-1}\right.\right.\\
&\left.\left.-\frac{1-b}{2b}\left(1-t\right)^{2}\sum_{n=2}^{\infty}n\left(n-1\right)b^{n}t^{n-2}\left(1+O\left(\frac{1}{\ln N}\right)\right)\right]\right\}\,dt,
\end{align*}
where 
\begin{equation}
b=\frac{p}{p+1}, \label{beta}
\end{equation}
and we have used that
\begin{align*}
\left(1-bt\right)^{-1}=\sum_{n=0}^{\infty}\left(b\,t\right)^{n},\,\,\,\left(1-bt\right)^{-2}=b^{-1}\sum_{n=1}^{\infty}nb^{n}\,t^{n-1},\\
\,\,\,\left(1-bt\right)^{-3}=2b^{-2}\sum_{n=2}^{\infty}n\left(n-1\right)b^{n}\,t^{n-2},
\end{align*}
since $\varepsilon \in (0,1)$, $b \in (0,1)$, and $t\in[-\varepsilon,\varepsilon]$. If we change the variables as $N^{t}=u\,\omega^{m-1}$, where $\omega:=\left(\ln N\right)^{-1}$, and apply the binomial theorem, after some careful computations we get
\begin{align*}
I_{2}(N)=\,\omega\int_{\omega^{1-m}\exp\left(-\varepsilon /\omega\right)}^{\omega^{1-m}\exp\left(\varepsilon /\omega\right)}
&\exp\left\{-\frac{\left(1-b\right)u}{\left(m-1\right)!}\left[\,1+\left(b-\left(m-1\right)\right)\omega\,\ln\left(u \omega^{m-1}\right)\right.\right.\\
&\left.\left.+ O\left(\omega\,\ln\left(u \omega^{m-1}\right)\right)^{2}\right]\right\}\\
&\times \exp\left\{-\frac{\omega\, u}{\left(m-1\right)!}\left[d_{1}+O\left(\omega\,\ln\left(u \omega^{m-1}\right)\right)\right]\right\}\frac{du}{u},
\end{align*}
where 
\begin{equation}
d_{1}=\left(1-b^{2}\right)\left(m-1\right)-\frac{1-b}{b}-3b^{2}\left(1-b\right).\label{d1}
\end{equation}
Notice that, $N\rightarrow \infty$ implies $\omega \rightarrow 0^+$. We claim that we can replace the upper limit in the above expression by $\infty$. Let us rewrite $I_{2}(N)$ as
\begin{equation}
I_{2}(N)=\omega \left(\int_{\omega^{1-m}\exp\left(-\varepsilon /\omega\right)}^{1/\sqrt{\omega}}
+ \int_{1/\sqrt{\omega}}^{\omega^{1-m}\exp\left(\varepsilon /\omega\right)}\right). \label{eint5}
\end{equation}
The second integral of (\ref{eint5}) is easily bounded by 
$O\left(\sqrt{\omega}\; e^{-(1-b)/\left(m-1\right)!\sqrt{\omega}}\right)$.
Let us denote $I_{21}(\omega)$ the first integral of (\ref{eint5}). We expand the exponentials and get
\begin{align*}
I_{21}(\omega)=\int_{\omega^{1-m}\exp\left(-\varepsilon /\omega\right)}^{1/\sqrt{\omega}}\frac{e^{-\left(1-b\right)u/\left(m-1\right)!}}{u}
&\left[1-\frac{1-b}{\left(m-1\right)!}\left(b-\left(m-1\right)\right)u\omega\,\ln\left(u\omega^{m-1}\right)\right.\nonumber\\
&\left.\,\,\,\,\,\,-\frac{d_{1}}{\left(m-1\right)!}\,u\,\omega\,\left(1+O\left(\omega\ln\left(u\,\omega^{m-1}\right)\right)\right)\right]du.
\end{align*}
We write the integral above as
\begin{equation}
I_{21}(\omega)=\int_{\omega^{1-m}\exp\left(-\varepsilon /\omega\right)}^{\infty}-\int_{1/\sqrt{\omega}}^{\infty}.\label{20}
\end{equation}
Again, the second integral of (\ref{20}) is easily bounded by $O\left(\sqrt{\omega}\,e^{-\left(1-b\right)/\left(m-1\right)!\sqrt{\omega}}\right)$ as $\omega\rightarrow 0^{+}$, and our claim is proved. It is now an easy exercise to evaluate $I_{2}(N)$. We have
\begin{align*}
I_{2}(N)=&\,\varepsilon+\left(m-1\right)\omega \ln \omega+\left[\,\ln\left(m-1\right)!+\ln\left(p+1\right)- \gamma\,\right] \omega\nonumber\\
&-\left(m-1\right)\left(b-\left(m-1\right)\right) \omega^{2}\ln \omega\nonumber\\
&+\left[\left(b-\left(m-1\right)\right)\left(\gamma-\ln\left(m-1\right)!-\ln\left(p+1\right)-d_{1}\left(1-b\right)\right)\right]\omega^{2}\nonumber\\
&+O\left(\omega^{3}\left(\ln \omega\right)^{2}\right),
\end{align*}
(where $b$ and $d_{1}$ as defined in (\ref{beta}) and (\ref{d1}) respectively). Notice that the error term in the above dominates the previously mentioned term $O\left(\sqrt{\omega}\,e^{-\left(1-b\right)/\left(m-1\right)!\sqrt{\omega}}\right)$ as $\omega\rightarrow 0^{+}$.\\
Now, we turn our attention to $I_{3}(N)$ of (\ref{I3}). As we will see the leading term is enough. The idea is that one can replace the integrand of (\ref{I3}) with $\left[\,-M_{m}(N;s)\,\right]$ and then by the quantity
\begin{equation*}
N_{m}(N;s):=\sum_{j=2}^N  \left[S_{m}\bigg(\frac{(\ln N)^{p+1}}{ \left(\ln j\right)^{p}}s\bigg) \exp{\bigg(-\frac{(\ln N)^{p+1}}{ \left(\ln j\right)^{p}}s\bigg)} \right].
\end{equation*}
For a rigorous approach see, \cite{DP}. Hence as $N\rightarrow \infty$
\begin{equation*}
I_{3}(N)=\int_{1+\varepsilon}^{\infty} N_{m}(N;s)\left[1+O\left(N_{m}(N;s)\right)\right] ds.
\end{equation*}
From the comparison of sums and integrals and Lemma 2.1 one easily arrives at
\begin{equation*}
I_{3}(N)=\sum_{k=0}^{m-1}\frac{\left(\ln N\right)^{k}}{k!}\,\int_{1+\varepsilon}^{\infty} \frac{s^{k}N^{1-s}}{1+ps}
\left[1+O\left(\ln N\right)\right]ds.
\end{equation*}
Substitute  $s=1-t$ and apply the Lapace method for integrals yields
\begin{equation}
I_{3}(N)=\frac{\left(1+\varepsilon\right)^{m-1}}{\left(1+p\right)\left(m-1\right)!\,\omega^{m-2}}\,e^{-\varepsilon/\omega}
\left[1+O\left(\frac{1}{\omega}\right)\right]
\label{I3RESULT}
\end{equation}
as $\omega \rightarrow 0^{+}$ and as we have set $\omega=\left(\ln N\right)^{-1}$. The reader now observes that the leading term of $I_{3}(N)$ is dominated by the sixth term of $I_{2}(N)$ as $N\rightarrow \infty$. We finish our approach by estimating the integral $I_{1}(N)$ of (\ref{I1}). For any given $\varepsilon \in (0,1)$ it is easy to see that
\begin{align*}
I_{1}(N)&= \int_0^{1 - \varepsilon} \exp \left[
\sum_{j=2}^{N}\ln \left[1-S_{m}\bigg(\frac{(\ln N)^{p+1}}{ \left(\ln j\right)^{p}}s\bigg) \exp{\bigg(-\frac{(\ln N)^{p+1}}{ \left(\ln j\right)^{p}}s\bigg)}\right] \right] ds\nonumber\\
&< \exp \left[-\sum_{k=0}^{m-1}\left[\frac{\left(1-\varepsilon\right)^{k}\ln N^{k\,(p+1)}}{\left(m-1\right)!}
\left( \sum_{j=2}^{N} \left(\ln j\right)^{-kp}e^{-\left(1 - \varepsilon \right)\frac{\left(\ln N\right)^{p+1}}{\left(\ln j\right)^{p}}}\right)\right]\right].
\end{align*}
From the comparison for sums and integrals it follows that (as $N\rightarrow \infty$)
\begin{equation*}
\sum_{j=2}^{N} \left(\ln j\right)^{-kp}e^{-\left(1 - \varepsilon \right)\frac{\left(\ln N\right)^{p+1}}{\left(\ln j\right)^{p}}}\sim
\int_{j=2}^{N} \left(\ln x\right)^{-kp}e^{-\left(1 - \varepsilon \right)\frac{\left(\ln N\right)^{p+1}}{\left(\ln x\right)^{p}}}dx.
\end{equation*}
Since $1-\varepsilon$ is strictly positive \textit{it is safe} to apply Lemma 2.1 and easily arrive at the inequality
\begin{align*}
I_{1}(N)&
< \exp\left[-\sum_{k=0}^{m-1}\frac{\left(1-\varepsilon\right)^{k}}{\left(1+p\left(1-\varepsilon\right)\right)\left(m-1\right)!}\,\frac{e^{\varepsilon/\omega}}{\omega^{k}}\,\left(1 + M_1\,\omega \right) \right]\nonumber\\
&=\exp\left[-\frac{1}{\left(1+p\left(1-\varepsilon\right)\right)\left(m-1\right)!}\,\frac{\omega^{m}-\left(1-\varepsilon\right)^{m}}{\omega^{m-1}\left(\omega-\left(1-\varepsilon\right)\right)}\,e^{\varepsilon/\omega}\,\left(1 + M_1\,\omega \right) \right],
\end{align*}
where $M_1$ is a positive constant. Since $\omega \rightarrow 0^{+}$ and $\varepsilon \in (0,1)$ we have
\begin{equation*}
I_{1}(N) << \frac{\left(1+\varepsilon\right)^{m-1}}{\left(1+p\right)\left(m-1\right)!\,\omega^{m-2}}\,e^{-\varepsilon/\omega},
\end{equation*}
for sufficiently large $N$, $m=1,2,3,\cdots.$ \\
Now \textbf{Main result I} follows immediately. It is notable that the \textit{third term} of $A_{N}=\sum_{j=2}^{N}\left(\ln j\right)^{-p}$ contributes to the average of $T_{m}(N)$.

\smallskip

\bigskip

\smallskip

\textbf{Proof of main result II.}
From (\ref{15}), (\ref{13}), and (\ref{SD4}) we have  
\begin{align}
&E[\,T_m(N)\left(T_m(N)+1\right)\,]=2\,N^{2}\bigg(\left(\ln N\right)^{2} +2 p\,\ln N+\left(p^{2}+2p\left(p+1\right)\right)+ O\left( \frac{1}{\ln N}\right)\bigg)\nonumber \\
&\,\,\,\,\,\,\,\,\,\,\,\,\,\,\,\,\,\times \int_{0}^{\infty}s\left\{1-\exp \Bigg( \sum_{j=2}^{N}\ln
\left[1-S_{m}\bigg(\frac{(\ln N)^{p+1}}{ \left(\ln j\right)^{p}}s\bigg) \exp{\bigg(-\frac{(\ln N)^{p+1}}{ \left(\ln j\right)^{p}}s\bigg)} \right]\Bigg)\right\} ds.\label{RII}
\end{align}
Let us denote $\tilde{Q}_{m}(N;\alpha)$ the integral above. Then, for any given $\varepsilon \in (0,1)$ we have
\begin{align*}
\tilde{Q}_{m}(N;\alpha)= \left[\,\frac{1}{2}+\varepsilon+\varepsilon^{2} -I_4 (N)-I_5 (N)+I_6 (N)\,\right],
\end{align*}
where
\begin{align}
I_4 (N):&= \int_0^{1-\varepsilon} s\,e^{M_{m}(N;s)}\, ds,\nonumber\\
I_{5}(N): &= \int_{1-\varepsilon}^{1+\varepsilon } s\,e^{M_{m}(N;s)}\, ds,\label{I5} \\
I_{6}(N): &=\int_{1+\varepsilon}^{\infty } s\left[1-e^{M_{m}(N;s)}\right]\, ds,\nonumber
\end{align}
and $M_{m}(N;s)$ is given in (\ref{AsN}).
If we treat $I_5(N)$ as we treated $I_2(N)$ and with a little patiences and paper, one finally arrives at
\begin{align*}
I_{5}(N)=&\varepsilon+\frac{\varepsilon^{2}}{2}+\left(m-1\right)\omega\ln \omega+\left[\ln\left(m-1\right)!+\ln\left(p+1\right)-\gamma\right]\omega
-\frac{\left(m-1\right)^{2}}{2} \omega^{2}\ln^{2} \omega\nonumber\\
&\,\,\,+\left(m-1\right)\left[\left(m-1\right)-\frac{p}{p+1}-\ln\left(m-1\right)!-\ln\left(p+1\right)+\gamma\right]\omega^{2}\ln \omega\nonumber\\
&\,\,\,+\left[\left(b-\left(m-1\right)\right)\left(\gamma-\ln\left(m-1\right)!-\ln\left(p+1\right)-d_{1}\left(1-b\right)\right)\right.\nonumber\\
&\left.\,\,\,\,\,\,\,\,\,-\frac{1}{2}\left(\gamma^{2}+\frac{\pi^{2}}{6}\right)+\gamma \left(\ln\left(m-1\right)!+\ln\left(p+1\right)\right)\right.\nonumber\\
&\left.\,\,\,\,\,\,\,\,\,+\frac{1}{2}\left(\ln\left(m-1\right)!+\ln\left(p+1\right)\right)^{2}\right]\omega^{2}+O\left(\omega^{3}\left(\ln\omega\right)^{2}\right),
\end{align*}
(where $b$ and $d_{1}$ as defined in (\ref{beta}) and (\ref{d1}) respectively). With similar steps as in Main result I one has that $I_{4}(N)$ and $I_{6}(N)$ are negligible compared to the \textit{eighth} of $I_{5}(N)$. Now Main result II follows immediately by invoking (\ref{RII}).

\smallskip

\bigskip

\smallskip

\textbf{Proof of main result III.} The proof follows immediately from the identity
\begin{equation*} 
V[\,T_m(N)\,]=E[\,T_m(N)\left(T_m(N)+1\right)\,]-E[\,T_m(N)\,]-E[\,T_m(N)\,]^{2}
\end{equation*}
by invoking Main results I and II.
\smallskip

\bigskip

\smallskip

\textbf{Proof of main result IV.} P. Neal \cite{N} has established a general theorem regarding the limit distribution of $T_m(N)$ (appropriately normalized) as $N \to \infty$,
where $\pi_N = \{p_{N1}, p_{N2},...,p_{NN} \}$, $N = 1, 2,...$, is a sequence of (sub)probability measures, not necessarily of the form (\ref{8}).

\smallskip

\textbf{Theorem N.} Suppose that there exist sequences $\{b_N\}$ and $\{k_N\}$ such that $k_N / b_N \rightarrow 0$
as $N \rightarrow \infty$ and that, for $y \in \mathbb{R}$,
\begin{equation}
\Lambda_N(y\,;m) := \frac{b_N^{m-1}}{\left(m-1\right)!} \sum_{j=1}^N p_{Nj}^{m-1}\exp\bigg(-p_{Nj} \left(b_N + y k_N\right) \bigg) \rightarrow g(y),
\quad
N \rightarrow \infty,
\label{N1}
\end{equation}
for a nonincreasing function $g(\cdot)$ with $g(y) \rightarrow \infty$ as $y \rightarrow -\infty$ and
$g(y) \rightarrow 0$ as $y \rightarrow \infty$. Then
\begin{equation}
\frac{T_{m}(N) - b_N}{k_N} \overset{D}{\longrightarrow} Y,
\qquad
N \rightarrow \infty,
\label{N2}
\end{equation}
where $Y$ has distribution function
\begin{equation}
F(y) = P\{ Y \leq y \} = e^{-g(y)},
\qquad
y \in \mathbb{R}.
\label{N222a}
\end{equation}.

Theorem N \textit{does not indicate at all} how to choose the sequences $\{b_N\}$ and $\{k_N\}$. Here our asymptotic formulas can help. In particular, we will choose
\begin{equation}
b_N = N \ln N + (m-1) N \ln\ln N
\qquad  \text{and} \qquad
k_N = N
\label{SD0}
\end{equation}
and  for all $y \in \mathbb{R}$ we will prove that
\begin{equation}
P\left\{\frac{T_m(N) - N \ln N - (m-1) N \ln\ln N}{N} \leq y \right\}\rightarrow \exp\left(-\frac{e^{-(y-p)}}{(p+1) (m-1)!}\right)
\label{SD1}
\end{equation}
as $N \rightarrow \infty$, which is equivalent to Main result IV. Under the choice of (\ref{SD0}), $\Lambda_N(y\,;m)$ of (\ref{N1}) satisfies, as $N \rightarrow \infty$,
\begin{equation}
\Lambda_N(y\,; m) \sim \frac{(N \ln N)^{m-1}}{(m-1)!}
\sum_{j=2}^N \left(\frac{a_j}{A_N}\right)^{m-1} e^{-(a_j / A_N) (N \ln N + (m-1) N \ln\ln N + N y)}
\label{SD2}
\end{equation}
where 
\begin{equation*}
a_j = \frac{1}{(\ln j)^p}\quad\text{and}\quad
A_N = \sum_{j=2}^N \frac{1}{(\ln j)^p} = \frac{N}{(\ln N)^p} + \frac{pN}{(\ln N)^{p+1}} + O\left( \frac{N}{(\ln N)^{p+2}} \right).
\end{equation*}
Hence, \eqref{SD2} yields
\begin{equation}
\Lambda_N(y\,; m) \sim \frac{(\ln N)^{(p+1)(m-1)}}{(m-1)!} \, S_N(y),
\label{SD5}
\end{equation}
where
\begin{equation}
S_N(y) := \sum_{j=2}^N \frac{1}{(\ln j)^{p(m-1)}} \exp\left(-\frac{(\ln N)^p (1 - p / \ln N) (\ln N + (m-1) \ln\ln N + y)}{(\ln j)^p}\right).
\label{SD6}
\end{equation}
Now,
\begin{equation}
S_N(y) \sim I_N(y)
\label{SD7}
\end{equation}
where
\begin{equation}
I_N(y) := \int_2^N \frac{1}{(\ln x)^{p(m-1)}} \exp\left(-\frac{(\ln N)^p (1 - p / \ln N) (\ln N + (m-1) \ln\ln N + y)}{(\ln x)^p}\right) dx.
\label{SD8}
\end{equation}
By substituting $u = \ln x$ in the above integral we get
\begin{equation}
I_N(y) := \int_2^M \frac{1}{u^{p(m-1)}} \exp\left(-\frac{B}{u^p} + u\right) du,
\label{SD9}
\end{equation}
where for typographical convenience we have set
\begin{equation}
B := \omega^{-\left(p+1\right)} \left(1 - p\,\omega\right) \left(1 - \frac{(m-1)\omega \ln \omega}{M} + y\,\omega\right)
\quad \text{and} \quad
\omega:= \left(\ln N\right)^{-1}
\label{SD10}
\end{equation}
so that $B \rightarrow \infty$ and $\omega \rightarrow 0^{+}$ as $N \rightarrow \infty$.

Next, in the integral of \eqref{SD9} we substitute $u = B^{1 / (p+1)} t$ and obtain
\begin{equation}
I_N(y) \sim B^{1 - \frac{pm}{p+1}} \int_0^{\theta} \frac{1}{t^{p(m-1)}} \, e^{B^{1 / (p+1)} \phi(t)} dt,
\label{SD11}
\end{equation}
where
\begin{equation}
\theta := \omega^{-1} B^{-1 / (p+1)}
\quad \text{and} \quad
\phi(t) := t - \frac{1}{t^p}.
\label{SD12}
\end{equation}
The integral in the right-hand side of \eqref{SD11} can be treated as a Laplace integral \cite{B-O}, where the large parameter is $B^{1 / (p+1)}$.
Since $\phi(t)$ is strictly increasing, the main contribution to the asymptotics of this integral comes from the endpoint $\theta$ (notice that
$\theta \sim 1$ as $N \to \infty$). Thus, by applying the standard analysis of Laplace integrals, after some straightforward algebraic manipulations \eqref{SD11} becomes
\begin{equation}
I_N(y) \sim M^{-(p+1)(m-1)} \frac{e^{-(y-p)}}{(p+1)}.
\label{SD13}
\end{equation}
Finally, by combining \eqref{SD13} with \eqref{SD10}, \eqref{SD7}, and \eqref{SD5} we obtain
\begin{equation}
\Lambda_N(y\,; m) \sim \frac{e^{-(y-p)}}{(p+1) (m-1)!}
\label{SD14}
\end{equation}
and the proof is finished by invoking Theorem $N$.

\end{document}